\documentclass[11pt,english]{article}
\usepackage{babel,amsfonts}
\usepackage{indentfirst,graphics,graphicx}
\usepackage[dvipsnames,usenames]{xcolor}
\usepackage{amsfonts,amssymb,latexsym,amsmath}
\usepackage{a4wide}
\usepackage[latin1]{inputenc}
\usepackage{cancel,soul} % para tachar fórmulas

\usepackage[margin=1in]{geometry}

\def\u{{\textbf {\textit u}}}
\def\d{{\textbf {\textit d}}}
\def\x{{\textbf {\textit x}}}

\def\W{{\textbf {\textit W}}}

\newtheorem{theorem}{Theorem}[section]

\newtheorem{corollary}[theorem]{Corollary}

\newtheorem{remark}{Remark}[section]

\begin{document}

\title{A uniqueness and regularity  criterion for Q-tensor models with Neumann boundary conditions.}

% If your title is too long to fit on the running head:
% \title[short title]{Title of your article}

% If your title is too long to fit on one line:
% \title[short title]{Title first line  \\ title second line}

%\author{Francisco Guill{\'e}n-Gonz{\'a}lez and Mar{\'\i}a {\'A}ngeles Rodr{\'\i}guez-Bellido}

\author{Francisco Guill{\'e}n-Gonz{\'a}lez
  and Mar{\'\i}a {\'A}ngeles Rodr{\'\i}guez-Bellido \\
  Dpto.~de Ecuaciones Diferenciales y An\'alisis Num\'erico and IMUS, \\
  Universidad de Sevilla \\
Aptdo.~de Correos 1160, 41080 Sevilla, Spain \\
  \texttt{guillen@us.es} and \texttt{angeles@us.es}}
  
\date{}
\maketitle
% If there are two authors:
% \author{First author and second author}

% If there are three or more authors:
% \author{First author, second author,...,nth author and last author}

% If the authors' names do not fit on one line:
% \author[abbreviated list]{first line \\ second line}
% In the optional argument put something like first author et al.

\begin{abstract}
We give a regularity criterion for a $Q$-tensor system 
%with stretching term 
modeling  a nematic Liquid Crystal, under homogeneous Neumann boundary conditions for the tensor $Q$. Starting of a criterion  only imposed on the velocity field $\u$  two results are proved; the uniqueness of weak solutions and the global in time weak regularity for the time derivative $(\partial_t \u,\partial_t Q)$.  This paper extends the work done in \cite{GGRBRM} for a nematic Liquid Crystal model formulated in  $(\u,\d)$, where $\d$ denotes the orientation  vector of the liquid crystal molecules.
\end{abstract}

{\bf Key words:} Nematic  liquid crystal, Q-tensor model, regularity criteria, uniqueness criteria, Neumann boundary conditions.

\let\thefootnote\relax\footnote{{\bf AMS 2010 subject classification}: 35B65, 35K51, 35Q35, 76A15, 76D03}

\let\thefootnote\relax\footnote{The authors have been partially supported by MINECO Project~MTM2012-32325 of Spain Government.}

%\newpage
%
% 
%
%{\bf Abbreviated version of the title for the running head:} \par 
%
%\begin{center}
%{\Large Uniqueness-regularity Q-tensor criterion}
%\end{center}
%
%\bigskip
%
%{\bf Corresponding author:} \par 
%
%\begin{center}
%Mar{\'\i}a {\'A}ngeles Rodr{\'\i}guez-Bellido \\
%  Dpto.~de Ecuaciones Diferenciales y An\'alisis Num\'erico, \\
%  Facultad de Matem{\'a}ticas, Universidad de Sevilla \\
%Aptdo.~de Correos 1160 \\
%41080 Sevilla, Spain \\
%\texttt{e-mail: angeles@us.es}
%\end{center}
%
% 
%
%\newpage

% ====================================================================
% Section 1
% ====================================================================

\section{Introduction}
\subsection{The model.}

Liquid crystals are intermediate phases of matter with properties from both solid and liquid states. The macroscopic properties come from the liquid behavior and are modeled by using the velocity and pressure $(\u,p)$. The microscopic structure enters into the model through the molecules of liquid crystals. These molecules influence the behavior of the matter and such an influence can be modeled by using different kinds of unknowns, which depend on different theories and the type of liquid crystal. 

In the nematic case, where the molecules are arranged by layers and every molecule is oriented equally, two main theories appear: the Oseen-Frank theory, where the microscopic structure is modeled by using a director vector $\d$, which is the average of the main orientation of the rod-like molecules of liquid crystals, and the Landau-De Gennes theory, where the vector $\d$ is replaced by a tensor $Q$.

 The tensor $Q$ is related to the second moment of a
probability measure $\mu(\x,\cdot): \, \mathcal{L}(\mathbb{S}^2)\rightarrow [0,1]$ for each $\x\in \Omega$, being $\mathcal{L}(\mathbb{S}^2)$ the family of Lebesgue
measurable sets on the unit sphere. For any $A\subset \mathbb{S}^2$, $\mu(\x,A)$ is the probability that
the molecules with centre of mass in a very small neighborhood of the point $\x \in \Omega$ are pointing in a direction
contained in $A$.
From a physical point of view, this probability (cf.~\cite{Gennes,MZ}) must satisfy $\mu(\x,A)=\mu(\x,-A)$ in order to reproduce the so-called
``head-to-tail'' symmetry. As a consequence, the first moment of the probability measure vanishes, that is
$$\langle p \rangle (\x)=\displaystyle\int_{\mathbb{S}^2} p_i \, d\mu(\x,p)=0.$$
Then, the main information on $\mu$ comes from the second moment tensor
$$M(\mu)_{ij}=\displaystyle\int_{\mathbb{S}^2} p_i \, p_j \,
d\mu(p),\quad i,j=1,2,3.$$
As a consequence, $M(\mu)=M(\mu)^t$ and $tr(M)=1$. If the orientation of the
molecules is equally distributed, then the distribution is
isotropic and $\mu=\mu_0$, $d\mu_0(p)=\frac{1}{4\pi} \, dA$
and $M(\mu_0)=\frac{1}{3} \, \mathbb{I}$.  The deviation of the
second moment tensor from its isotropic value is therefore measured as:
$$
%\begin{equation}\label{d1}
Q=M(\mu)-M(\mu_0)=\displaystyle\int_{\mathbb{S}^2}
\left(p\otimes p - \displaystyle\frac{1}{3}\,
\mathbb{I}\right) \, d\mu(p),
$$
%\end{equation}
which is the definition for the tensor $Q$, and from which the symmetry and traceless for $Q$ is deduced. 

By following the Landau-De Gennes theory, a model to study the behavior of nematic liquid crystals filling  a bounded domain $\Omega\subset \mathbb{R}^{3}$, with  boundary $\partial\Omega$, is given by:
\begin{equation} \label{modelo-u-1}
\left\{\begin{array}{l}
D_t \u  -\nu \Delta \u
+\nabla p  =  \nabla \cdot \tau(Q) + \nabla \cdot \sigma (H,Q), \quad  \nabla \cdot \u  =0 \quad \mbox{in $\Omega \times
(0,T),$}\\
D_t Q -S(\nabla \u , Q)  = -\gamma\, {H}(Q)
\quad \mbox{in $\Omega \times
(0,T),$}
\end{array}\right.
\end{equation}
where $\u:(0,T)\times \Omega \rightarrow \mathbb{R}^3 $ is the velocity field, $p:(0,T)\times \Omega \rightarrow \mathbb{R}$ is the pressure and $Q:(0,T)\times \Omega \rightarrow  \mathbb{R}^{3\times 3}$ is a symmetric and traceless tensor. The operator
$D_t(\cdot)=\partial_t (\cdot) + (\u  \cdot \nabla)
( \cdot  ) $ is the material derivative, $\nu>0$ is the viscosity
coefficient, and $\gamma>0$ is a material-dependent elastic constant. The tensors $\tau=\tau(Q)$,
$\sigma=\sigma(H,Q)$ are defined by:
\begin{equation}\label{modelo-2}
\left\{\begin{array}{rcl}
\tau_{ij}(Q) & = & -\varepsilon
\, \left( \partial_j Q :
\partial_i Q  \right)= -\varepsilon \, \partial_j Q_{kl} \, \partial_i Q_{kl}  \quad (\varepsilon >0), \\
\sigma (H,Q) & = &  {H} \, Q-Q \, {H}   ,
\end{array}\right.
\end{equation}
where $H=H (Q)= - \varepsilon \, \Delta Q  +f(Q)$ and
\begin{equation}\label{fz}
f(Q)=a \, Q -b \, (Q^2-\dfrac{1}{3} \, tr(Q^2)) + c \vert Q \vert^2 \, Q \quad
\mbox{ with  $a$, $b\in \mathbb{R}$ and $c
>0$},
\end{equation}
hereafter, $\vert Q \vert^2 = Q : Q =\displaystyle\sum_{i,j=1}^3Q_{ij} \, Q_{ij}$ denotes the tensor euclidean norm.
The stretching term $S(\nabla \u  , Q) $ is given by
\begin{equation}\label{S-2}
S(\nabla \u  , Q) =   \W  \, Q^t-Q^t
\, \W
\end{equation}
where $\W=\frac{1}{2} \, (\nabla \u - (\nabla \u)^t)$ is the antisymmetric part of the gradient of $\u$. 

In \cite{Z1,Abels} the problem is solved in the space of traceless and symmetric matrices (in the whole $\mathbb{R}^3$ and in a bounded domain $\Omega \subset \mathbb{R}^3$, respectively). In \cite{Z1}, the existence of global in time weak solutions in $3D$, and strong regularity and weak-strong uniqueness results in $2D$ are proved. A generalized model considering more stretching effects is studied in \cite{Z2}. The analysis of a model considering a more general expression for $H(Q)$ (modifying the  $\Delta Q$-term  of $H(Q)$) can be seen in \cite{HD}.

However, a similar study can be made for a system that generalizes \eqref{modelo-u-1}-\eqref{S-2}, in which 
the definition of $S(\nabla \u  , Q) $ is replaced by:
\begin{equation}\label{def-S}
S(\nabla \u  , Q)  =   \nabla \u  \, Q^t-Q^t
\, \nabla \u 
\end{equation}
and the function $f(Q)$ in (\ref{fz}) is replaced by:
\begin{equation}\label{ff}
f(Q)=a \, Q - \displaystyle\frac{b}{3} \, \left(
Q^2+QQ^t+Q^tQ \right )+c  \, \vert Q \vert^2 \, Q \quad
\mbox{ with  $a$, $b\in \mathbb{R}$ and $c
>0$}.
\end{equation}
Although this model does not have the restrictions of symmetry and tracelessness, it can be studied  in the same way. In fact, the existence of global weak solutions and the weak/strong uniqueness of this model is proved in \cite{GGRBinprep3}, and the local regularity and uniqueness is obtained in \cite{GGRBinprep2}. This study tell us that these (physical) restrictions are not essential in the mathematical analysis. 
On the other hand, if  $S(\nabla \u,Q)$  is taken as in (\ref{S-2}), then any weak solution provides a symmetric tensor $Q$; meanwhile if $f(Q)$ is taken as in (\ref{fz}), any weak solution provides a  traceless  tensor $Q$ (\cite{GGRBinprep3}).

The PDE system is enclosed with the following initial and boundary conditions:
\begin{equation}\label{ini10}
\u \vert_{t=0} = \u _0, \quad  Q \vert_{t=0}=Q_0 \quad
\mbox{in $\Omega,$} 
\end{equation}
\begin{equation}\label{ini1}
\u  \vert_{\partial\Omega}={\bf 0}, \quad   \partial_{{\bf n}}Q
\vert_{\partial\Omega}=0  \quad  \mbox{in $(0,T)$.}
\end{equation}
The compatibility condition $\partial_{{\bf n}}Q_0
\vert_{\partial\Omega}=0$ must be satisfied.

In the present work, we focus on the problem without restrictions of symmetry and tracelessness (\ref{modelo-u-1})-(\ref{modelo-2}), (\ref{def-S})-(\ref{ini1}).   At the end of the work, the results proven for this model are extended to the model with tracelessness and symmetry (\ref{modelo-u-1})-(\ref{S-2}), (\ref{ini10})-(\ref{ini1}).

 The manuscript is organized as follows: Next subsection summarizes the results on regularity and uniqueness  for the problem (\ref{modelo-u-1})-(\ref{modelo-2}), (\ref{def-S})-(\ref{ini1}). In Subsection \ref{new-result}, the new results in this paper are stated: Theorems \ref{t1} and \ref{t2}, and Corollaries \ref{c1} and \ref{c2}. Section \ref{nh} is devoted to the proof of Theorem \ref{t1} where the spatial treatment of the boundary term is made in Subsection \ref{bt}. The proof of Theorem \ref{t2} is made in Section \ref{sec-reg}. The last section analyzes how to extend the results in Subsection \ref{new-result} to the problems imposing symmetry and/or tracelessness.

\subsection{Known  results on regularity and uniqueness}
The weak solution for these models is defined extending the Navier-Stokes case. Specifically,  $({\u},Q)$ is said a weak solution in $(0,T)$ 
of (\ref{modelo-u-1})  if  $$
\u  \in L^{\infty}(0,T;{\bf L}^2(\Omega)) \cap L^2(0,T;{\bf H}^1(\Omega)),
\quad
Q \in L^{\infty}(0,T;\mathbb{H}^1(\Omega)) \cap
L^2(0,T;\mathbb{H}^2(\Omega)),
$$
satisfying the $\u$-system $(\ref{modelo-u-1})_1$ in a variational setting and the $Q$-system $(\ref{modelo-u-1})_2$ point-wisely.
The existence of weak solutions is proved in \cite{Z1,Z2} for the whole $\mathbb{R}^3$ and in \cite{Abels,GGRBinprep2} for a bounded domain of $\mathbb{R}^3$,  considering  Dirichlet boundary condition for $\u$ together with  Dirichlet or  Neumann boundary conditions for $Q$, or periodic boundary conditions for $(\u,Q)$. These results are based on the following \textsl{energy equality} (see \cite{GGRBinprep2}):
$$
%\begin{equation}\label{eq3}
\displaystyle\frac{d}{dt} \left(
\displaystyle\frac{1}{2} \Vert {\u}
\Vert_{{\bf L}^2(\Omega)}^2 + \displaystyle\int_{\Omega}
\mathcal{E}(Q) \, d{\bf x} \right) + \nu \Vert \nabla {\bf
u} \Vert_{{\bf L}^2(\Omega)}^2 + \gamma \Vert {H}
\Vert_{\mathbb{L}^2(\Omega)}^2 = 0 ,
$$
where the free energy is  $\mathcal{E}(Q)=\displaystyle\int_\Omega\left(\frac{\varepsilon}{2} \vert
\nabla Q \vert^2 + F(Q)\right)$ with the functional $F(Q)$ defined
as
\begin{equation}\label{Fgrande}
F(Q)=\displaystyle\frac{a}{2} \, \vert Q \vert^2
-\displaystyle\frac{b}{3} \, (Q^2 : Q)
+\displaystyle\frac{c}{4} \, \vert Q \vert^4 .
\end{equation}
 Note that,  in the case of $f(Q)$ defined by (\ref{ff}), it is easy to check that $F'(Q)=f(Q)$ and the tensor $H$ is  the variational derivative in $L^2(\Omega)$ of  $\mathcal{E}(Q)$, that is $H=\displaystyle\frac{\delta \mathcal{E}(Q)}{\delta Q}$ (see \cite{GGRBinprep2}).

%\end{equation}
%This energy equality has a force term depending on $\partial_t Q_{\partial\Omega}$ in the case of Dirichlet time-dependent boundary conditions for $Q$. 
On the other hand, in the case of bounded domains, the $Q$-system  $(\ref{modelo-u-1})_2$  satisfies a maximum principle  (see \cite{GGRBinprep2}), hence in particular 
\begin{equation}\label{max}
Q\in L^{\infty}(0,T;\mathbb{L}^\infty(\Omega)).
\end{equation}
 However, the uniqueness of weak solutions needs additional regularity for $\nabla \u$ and $\Delta Q$, which corresponds to the criterion proved by Berselli in \cite{Be} for the Navier-Stokes system. Concretely, the following result is proved in \cite{GGRBinprep2}:
\begin{theorem}[{\bf Uniqueness criteria}]\label{uniqueness-t}
Assume  $({\u}_0,Q_0)\in {\bf L}^2(\Omega) \times
\mathbb{H}^1(\Omega)$. Let $({\u},Q)$ be a weak solution in $(0,T)$ 
of problem (\ref{modelo-u-1})   such that $\nabla {\u}$ and $\Delta Q$ have the additional regularity
\begin{equation}\label{ar-u}
  \nabla {\u} \in
L^{2q/(2q-3)}(0,T;{\bf L}^q(\Omega))  \quad \mbox{for } 2\le q\le 3,
\end{equation} 
\begin{equation}\label{ar-Q}
  \Delta Q \in L^{2s/(2s-3)}(0,T;\mathbb{L}^s(\Omega)) \quad \mbox{for } 2\le s\le 3.
\end{equation}
Then, this solution coincides in
$(0,T)$ with any weak solution associated to the same data.
\end{theorem}
We can study two types of regularity for weak solutions $(\u,Q)$ of (\ref{modelo-u-1}):
\begin{itemize}
\item Strong regularity (as in the Navier-Stokes framework):
 \[
(\mbox{St})  \left\{\begin{array}{ll}
\u  \in L^{\infty}(0,T;{\bf H}^1(\Omega)) \cap L^2(0,T;{\bf H}^2(\Omega)),
& \partial_t \u  \in L^2(0,T;{\bf L}^2(\Omega)), \\
\noalign{\vspace{-2ex}}\\
Q \in L^{\infty}(0,T;\mathbb{H}^2(\Omega)) \cap
L^2(0,T;\mathbb{H}^3(\Omega)), &
\partial_t Q \in L^{\infty}(0,T;\mathbb{L}^2(\Omega)) \cap
L^2(0,T;\mathbb{H}^1(\Omega)).
\end{array}\right.
\]

\item Weak regularity for $(\partial_t\u ,\partial_tQ)$ or ``weak-t" solution:
\[
\hbox{(\mbox{w-t})}  \left\{\begin{array}{ll}
\partial_t \u  \in L^{\infty}(0,T;{\bf L}^2(\Omega)) \cap L^2(0,T;{\bf H}^1(\Omega)), &
\u  \in L^{\infty}(0,T;{\bf H}^1(\Omega)),
\\
\noalign{\vspace{-2ex}}\\
\partial_t Q \in L^{\infty}(0,T;\mathbb{H}^1(\Omega)) \cap
L^2(0,T;\mathbb{H}^2(\Omega)),
& Q \in L^{\infty}(0,T;\mathbb{H}^2(\Omega)).
\end{array}\right.
\]

\end{itemize}
%Both types of regularity coincides in the Navier-Stokes case and, in general, if no stretching terms are considered in the problem (\ref{modelo-u-1}).
Under homogeneous Neumann or Dirichlet boundary conditions for $Q$,  the existence and uniqueness of global in time strong solutions are obtained in \cite{GGRBinprep3} if the viscosity $\nu$ is large enough. Otherwise,  when $Q\vert_{\partial\Omega}= 0$  is imposed, and  the stretching term $S(\cdot,\cdot)$ is defined as in  \eqref{S-2}, then the existence and uniqueness of strong solutions are obtained in \cite{GGRBinprep3}; either local in time,  or global one if $\nabla \u$ satisfies the regularity criterion given below in (\ref{reg-ad-u}) (but regularity criterion for $\Delta Q$ given in (\ref{reg-ad-Q}) can be avoided). The existence of local in time strong solutions when $Q\vert_{\partial\Omega}= 0$ are also obtained via a fixed-point argument in \cite{Abels2}.

\begin{remark}
When the model (\ref{modelo-u-1})-(\ref{ini10}) is considered in the whole space $\Omega=\mathbb{R}^3$, some regularity criteria to obtain global in time strong solution (in the sense of $(St)$) are given in \cite{FO}. One of them imposes that $\nabla \u$ has the regularity appearing in (\ref{ar-u}) for $2<q\le 3$. 
Following Remark 3.2 in \cite{GGRBinprep3}, if we consider the model (\ref{modelo-u-1})-(\ref{ff}) with space-periodic boundary conditions for $(\u,Q)$, the regularity criterion only given for $\nabla \u$ in (\ref{ar-u}) implies the global in time strong solution in the sense of $(St)$ (no hypothesis for the symmetry of $S$ is assumed). Observe that in \cite{GGRBinprep3} the case $q=2$ in (\ref{ar-u})  is also considered.
\end{remark}

 However, as far as we know, the previous strong regularity results when $Q\vert_{\partial\Omega}= 0$ cannot be extended to  either non-homogeneous Dirichlet or  Neumann boundary conditions for $Q$. In these cases, some boundary integrals  appearing in the  strong estimates argument
 %for $\u$ and $Q$ in $L^{\infty}(0,T;{\bf H}^1(\Omega))\cap L^2(0,T;{\bf H}^2(\Omega))$ and $L^{\infty}(0,T;{\bf H}^2(\Omega))\cap L^2(0,T;{\bf H}^3(\Omega))$, respectively, 
 do not vanish and it is not clear  how to bound them. 
 %The study for can be deduced from the study given in \cite{SL,JBC} for a nematic model. 

 In order to circumvent this difficulty, the weak-t concept (w-t) is considered in \cite{Abels,GGRBinprep3}, obtaining weak-t solutions either local in time for any data or  global in time assuming some regularity criteria \cite{GGRBinprep3}, as can be summarized in the following result:

\begin{theorem}[{\bf Regularity criteria for global in time weak-t solution}]\label{rc2}   Assume $(\u _0,Q_0)\in {\bf H}^2(\Omega)\times \mathbb{H}^3(\Omega)$ and let $(\u ,Q)$ be a weak solution  in $(0,T)$ of problem  (\ref{modelo-u-1})
having the additional regularity:
\begin{equation}\label{reg-ad-u}
\nabla \u  \in L^{2q/(2q-3)}(0,T;{\bf L}^q(\Omega)),  \quad 3/2 \le q \le 3,
\end{equation}
\begin{equation}\label{reg-ad-Q}
 \Delta Q \in L^{2s/(2s-3)}(0,T;\mathbb{L}^s(\Omega)), \quad 3/2 \le  s \le 3.
\end{equation}
Then,   this weak solution $(\u ,Q)$ is a weak-t solution in
$(0,T)$ and coincides with any weak solution in
$(0,T)$ associated to the same data.
\end{theorem}
%
%Theorem \ref{rc2} can be extended to the case of 
%time-dependent Dirichlet boundary conditions for $Q$, that is $
% Q\vert_{\Gamma}=Q_{\Gamma}$ with $Q_{\Gamma} = Q_{\Gamma}(t)$ (see \cite{GGRBinprep3}).

For some models of nematic liquid crystals without stretching terms (see \cite{GGRBRM}), the criteria (\ref{ar-u}) and (\ref{ar-Q}) to obtain uniqueness of weak solutions or (\ref{reg-ad-u}) and (\ref{reg-ad-Q}) to obtain strong regularity could be replaced by the following criteria of Serrin's type (\cite{Serrin}):
%\begin{equation}\label{new1}
\begin{equation} \label{reg-u-nablad}
 {\u} \in L^{2p/(p-3)}(0,T;{\bf
L}^p(\Omega))  \qquad \mbox{and}
 \qquad \nabla Q \in
L^{2r/(r-3)}(0,T;\mathbb{L}^r(\Omega)), \quad 3 \le p, \, r \le
+\infty.
\end{equation}
%\end{equation}
Moreover, when periodic boundary conditions for $(\u,Q)$ are considered, the previous regularity criteria  imposed for $\Delta Q$ or $\nabla Q$ are not necessary.
%But, as long as we know, the presence of the stretching term does not allow to replace either hypothesis on $\nabla \u$ by  $\u$, or  hypothesis on $\Delta Q$ by  $\nabla Q$. Moreover, other boundary conditions for $Q$ different to periodic ones introduce a boundary integral term which is not clear how to control.

\begin{remark}
When $Q\equiv0$  the Q-tensor model reduces to the classical  Navier-Stokes problem. In this case, the regularity hypotheses for $\u$ (\ref{reg-ad-u}) and $(\ref{reg-u-nablad})_1$ correspond to Berselli's and Serrin's regularity criteria for the Navier-Stokes equations (see \cite{Be} and \cite{Serrin}, respectively). Indeed, the upper constraint $q\le 3$ in (\ref{reg-ad-u}) comes from the terms depending on $Q$ (see \cite{GGRBinprep2,GGRBinprep3}), and therefore they do not appear if $Q\equiv 0$.
\end{remark}

\bigskip

In this work, under  homogeneous Neumann boundary conditions for $Q$, a non-hilbertian regularity for $\nabla Q$ is  obtained only imposing a regularity criteria  for  $\nabla \u$. This non-hilbertian regularity argument for $\nabla Q$ follows the steps given in \cite{GGRBRM} for nematic liquid crystals without stretching terms and periodic boundary conditions. Afterwards, we use this new regularity for $\nabla Q$ (and the regularity criterion for $\nabla \u$) to deduce the additional regularity for $\Delta Q$ considered in (\ref{ar-Q}) and (\ref{reg-ad-Q}) for the particular index $s=5/2$.  This  regularity for $\Delta Q$ is based on  a  non-hilbertian regularity result for a  heat-Neumann problem (see \cite{Sol64}) and the $L^{\infty}$-regularity for $Q$ given by a maximum principle result (\cite{GGRBinprep2}).
% are crucial and we do not how to extend it to another boundary conditions, different to periodic ones.

\subsection{The new results of this paper}\label{new-result}
\begin{theorem}\label{t1}
Let $(\u,Q)$ be a weak solution of  (\ref{modelo-u-1})-(\ref{modelo-2}), (\ref{def-S})-(\ref{ini1})  such that $\nabla \u$ satisfies (\ref{ar-u}) and $Q$ satisfies (\ref{max}). Then, $\nabla Q \in L^{\infty}(0,T;\mathbb{L}^3(\Omega)) \cap L^3(0,T;\mathbb{L}^{9}(\Omega)).$
\end{theorem}
The proof is given in Section \ref{nh}. 
\begin{theorem}\label{t2}
Let $\Omega \subset \mathbb{R}^3$ be a bounded domain with boundary $\partial\Omega$ of class $C^{2+\epsilon}$ for some $\epsilon>0$. 
Let $(\u,Q)$ be a weak solution of  (\ref{modelo-u-1})-(\ref{modelo-2}), (\ref{def-S})-(\ref{ini1}) with $Q$ satisfying (\ref{max}).
Assume $({\u}_0,Q_0)\in {\bf L}^2(\Omega) \times
\mathbb{B}_{\gamma}^{2-2/\gamma}(\Omega)$ and hypotesis (\ref{reg-ad-u}) for $\nabla \u$, where $\gamma =\min\left\{q,{2q}/({2q-3})  \right\} $ and  $3/2\le q\le 3$ is the exponent given in (\ref{reg-ad-u}). Moreover, under the following hypothesis for the tensor:
 $$\nabla Q \in L^{\infty}(0,T; \mathbb{L}^3(\Omega)),$$
 the solution $(\u,Q)$  satisfies the additional regularity
\begin{equation}\label{gamma}
Q \in L^{\gamma}(0,T;{\bf W}^{2,\gamma}(\Omega)) \quad\hbox{and}\quad \partial_t Q \in L^{\gamma}(0,T;{\bf L}^{\gamma}(\Omega)).
\end{equation}
\end{theorem}
Note that, since $3/2\le q\le 3$, then $3/2\le \gamma\le 5/2$ and $2/3 \le 2-2/\gamma\le 6/5$. 

\begin{remark} 
Here,  $B_{\gamma}^{2-2/\gamma}(\Omega)$ is a Besov space.  In fact, Besov spaces   $B_{\gamma}^{l}(\Omega)$ coincides with  the well-known Sobolev spaces $W^{l, \gamma}(\Omega)$ when the  derivability exponent $l$ is  non-integer (\cite{Sol64}).
  In order to define  a Besov space $B_{\gamma}^{l}(\Omega)$ for  integer $l$, it suffices to start with  a one-variable function $f:\mathbb{R}\to \mathbb{R}$ where the norm in $B_{\gamma}^l(\mathbb{R})$ is given by
$$
\Vert f \Vert_{B_{\gamma}^l(\mathbb{R})} = \Vert f \Vert_{L^{\gamma}(\mathbb{R})} + \left(
\displaystyle\int_{-\infty}^{\infty}
\displaystyle\int_0^{\infty}
\left\vert
\displaystyle\frac{\partial^{l-1}f(x+2h)}{\partial x^{l-1}} - 2 \, \displaystyle\frac{\partial^{l-1}f(x+h)}{\partial x^{l-1}}
+\displaystyle\frac{\partial^{l-1}f(x)}{\partial x^{l-1}}
\right\vert^{\gamma} \, \displaystyle\frac{1}{h^{1+{\gamma}}}
dh \, 
dx
\right)^{1/{\gamma}}. 
$$
Then, the generalization to $B_{\gamma}^{l}(\Omega)$ for a $n$-dimensional domain $\Omega$ can be made via a partition of unity argument associated to overlapping subsets of $\Omega$, more details can be found in \cite{Sol64}. 

 The choice of the initial data $Q_0$ in the Besov space $B_{\gamma}^{2-2/\gamma}(\Omega)$ responds to the use of the result given below in  Theorem~\ref{t-3}. 
\end{remark}

We give the proof of Theorem \ref{t2} in Section \ref{sec-reg}. 
\medskip 

Theorems~\ref{t1} and \ref{t2} imply, in particular,  that 
\begin{equation}\label{mars-4}
\Delta Q \in L^{\gamma}(0,T;\mathbb{L}^{\gamma}(\Omega))\quad \mbox{ for $3/2\le \gamma\le 5/2$}.
\end{equation}
Regularity appearing in (\ref{mars-4}) can be seen as the regularity criteria of (\ref{ar-Q}) or (\ref{reg-ad-Q}) if and only if $\gamma=5/2$, because this is the unique case in (\ref{ar-Q}) or (\ref{reg-ad-Q}) where $s=2s/(2s-3)=5/2$. Therefore, since in Theorem~\ref{t2}  $\gamma =\min\left\{q,{2q}/({2q-3})  \right\} $, then  $\gamma=5/2$ corresponds with $q=5/2$ in the regularity criterion for $\nabla \u$ given in (\ref{ar-u}) or (\ref{reg-ad-u}), that is   $\nabla \u \in L^{5/2}(0,T;{\bf L}^{5/2}(\Omega))$.
As a consequence, by using $q=\gamma=5/2$ the following corollaries are deduced,  improving Theorems \ref{uniqueness-t} and \ref{rc2}, respectively:
\begin{corollary}[Uniqueness criteria]\label{c1}
Assume  $({\u}_0,Q_0)\in {\bf L}^2(\Omega) \times
\mathbb{W}^{2-2/(5/2),5/2}(\Omega)$. Let $({\u},Q)$ be a weak solution
of   (\ref{modelo-u-1})-(\ref{modelo-2}), (\ref{def-S})-(\ref{ini1})  such that $\nabla
{\u}$ satisfies the regularity criterion:
\begin{equation}\label{ar-b}
  \nabla {\u} \in
L^{5/2}(0,T;{\bf L}^{5/2}(\Omega)).
\end{equation} Then, this solution coincides in
$(0,T)$ with any weak solution associated to the same data.
\end{corollary}

\begin{corollary}[Regularity criteria]\label{c2}
Let $(\u ,Q)$ be a weak solution  in $(0,T)$ of  (\ref{modelo-u-1})-(\ref{modelo-2}), (\ref{def-S})-(\ref{ini1}). If $(\u _0,Q_0)\in {\bf H}^2(\Omega)\times \mathbb{H}^3(\Omega)$ and $\nabla
{\u}$ has the additional regularity (\ref{ar-b}), 
%:
%\begin{equation}\label{reg-ad-2}
%\nabla \u  \in L^{5/2}(0,T;{\bf L}^{5/2}(\Omega)) .
%\end{equation}
then $(\u ,Q)$
is the unique weak-t solution of  (\ref{modelo-u-1})-(\ref{modelo-2}), (\ref{def-S})-(\ref{ini1})  in $(0,T)$.
\end{corollary}

\section{Proof of Theorem \ref{t1}.}\label{nh}
We can argue as in  \cite{GGRBRM} for a Nematic Liquid Crystal model (without stretching term and periodic boundary conditions). Indeed,  by taking $-\nabla \cdot \left(
\vert \nabla Q \vert^{p-2} \nabla Q
\right)$, $p\ge 2$, as test function in the $Q$-system of (\ref{modelo-u-1}), using the incompressibility equation $\nabla\cdot\u=0$ and the non-slip boundary condition ${\u}|_{\partial\Omega}=0$, we obtain:
\[
\begin{array}{l}
-\displaystyle\int_{\Omega} \left(\partial_t Q +(\u \cdot \nabla ) Q -\gamma \, \varepsilon \Delta Q \right): [
\nabla \cdot ( \vert \nabla Q \vert^{p-2} \, \nabla
Q )] d\x \\
\noalign{\vspace{-1ex}}\\
\qquad 
 =
\displaystyle\frac{1}{p} \displaystyle\frac{d}{dt} \Vert \nabla
Q \Vert_{\mathbb{L}^p(\Omega)}^p +\gamma \,\varepsilon \displaystyle\int_{\Omega} \vert \nabla Q \vert^{p-2} \vert D^2 Q \vert^2
 \,d\x
\\
\noalign{\vspace{-1ex}}\\
\qquad +\displaystyle\int_{\Omega} \, \vert \nabla Q \vert^{p-2} \, (\nabla {\u} \cdot \nabla)
Q  : \nabla Q \, d \x 
+ \displaystyle\frac{1}{p} \,
\displaystyle\int_{\Omega} ({\u} \cdot \nabla) (\vert
\nabla Q \vert^p) \, d \x\\
\noalign{\vspace{-1ex}}\\
\qquad+\gamma \,\varepsilon \left(\displaystyle\frac{2}{p}\right)^2 \, (p-2) \,
\displaystyle\int_{\Omega} \left\vert \nabla \left( \vert \nabla
Q \vert^{p/2} \right) \right\vert^2 \, d\x - \gamma \,\varepsilon 
\displaystyle\int_{\partial\Omega}
\vert \nabla Q \vert^{p-2} \,
(\partial_k(\nabla Q)  \, n_k ) :  \nabla Q
\, d\sigma_{{\small \x}}
\end{array}
\]
To treat the potential term $F'(Q)=f(Q)$, we use the argument given in \cite{GGRBinprep2} splitting $f(Q)=F_c'(Q)+F_e'(Q)$, 
where 
$$
F_c'(Q)=c  \, \vert Q \vert^2 \, Q
$$
is the derivate of the  convex part $F_c(Q)=\displaystyle\frac{c}{4} \, \vert Q \vert^4$ (see (\ref{Fgrande})) and 
$$
F_e'(Q)=a \,  Q  - \displaystyle\frac{b}{3}  \, (Q^2+ Q \, Q^t+Q^t \, Q)
$$ 
is the derivate of the rest of $F(Q)$. Thus, 
$$
\begin{array}{rcl}
-\displaystyle\int_{\Omega}
F_c'(Q) : \left[ \nabla \cdot \left(
\vert \nabla Q \vert^{p-2} \nabla Q \right) \right]
\, d\x  & = & 
\displaystyle\int_{\Omega} \vert \nabla Q \vert^{p-2}
\nabla \left( F_c'(Q) \right) :
\nabla Q \, d\x\\
\noalign{\vspace{-1ex}}\\
& = & c \,   \displaystyle\int_{\Omega}
 \vert \nabla Q \vert^{p-2} \, \left(
\vert Q \vert^2 \vert \nabla Q \vert^2 + 2 \vert (Q:  \nabla Q )\vert^2
\right) \, d\x \ge 0,
\end{array}
$$
\[
\begin{array}{l}
-\displaystyle\int_{\Omega} F_e'(Q) :  \left[ \nabla \cdot ( \vert \nabla Q \vert^{p-2}
\nabla Q ) \right] d\Omega   =   \displaystyle\int_{\Omega}
\vert \nabla Q \vert^{p-2}
\nabla \left(F_e'(Q) \right) :
 \nabla Q \, d\x\\
 \noalign{\vspace{-1ex}}\\
\qquad  \le  a \, \Vert \nabla Q \Vert_{\mathbb{L}^p(\Omega)}^p +  2 \, b  \displaystyle\int_{\Omega} \vert \nabla Q \vert^p Q \, d\x\\
\noalign{\vspace{-1ex}}\\
\qquad \le   a \, \Vert \nabla Q \Vert_{\mathbb{L}^p(\Omega)}^p + \displaystyle\frac{c}{2}  
\displaystyle\int_{\Omega} \left[
\vert Q \vert^2 \vert \nabla Q \vert^p + 2 \, \vert \nabla Q \vert^{p-2} \, |(Q: \nabla Q)|^2
\right] \, d\x + C_{b,c} \, \Vert \nabla Q\Vert_{\mathbb{L}^p(\Omega)}^p
\end{array}
\]
(contrary to~\cite{GGRBRM}, the term depending on $F_e'(Q)$ is bounded without using a $L^\infty$-norm given by a maximum principle).
Therefore, considering homogeneous Neumann boundary conditions for $Q$, we obtain:
\[
\begin{array}{l}
\displaystyle\frac{1}{p}\displaystyle\frac{d}{dt} \Vert \nabla Q \Vert_{L^p(\Omega)}^p
  +   \gamma \, \varepsilon \displaystyle\int_{\Omega}
\vert D^2 Q \vert^2 \vert \nabla Q \vert^{p-2}
\, d\x + \gamma \,  \varepsilon \, (p-2) \, \left( \displaystyle\frac{2}{p}\right)^2 \Vert \nabla (\vert \nabla Q \vert^{p/2}) \Vert_{L^2(\Omega)}^2\\
\noalign{\vspace{-1ex}}\\
\qquad +  \displaystyle\frac{c \, \gamma}{2} \, \displaystyle\int_{\Omega} \left[
\vert Q \vert^2 \vert \nabla Q \vert^p + 2 \, \vert \nabla Q \vert^{p-2} \, |Q: \nabla Q|^2
\right] \, d\x \\
\noalign{\vspace{-1ex}}\\
\qquad \le   C \displaystyle\int_{\Omega} \vert \nabla {\u} \vert \vert \nabla Q \vert^p
\, d\x  -\displaystyle\int_{\Omega} S(\nabla {\u},Q):\nabla \cdot [\vert \nabla Q \vert^{p-2} \nabla Q] \, d\x\\
\noalign{\vspace{-1ex}}\\
\qquad +    \gamma \, a \, \Vert \nabla Q \Vert_{\mathbb{L}^p(\Omega)}^p + \gamma \, C_{b,c} \, \Vert \nabla Q\Vert_{\mathbb{L}^p(\Omega)}^p
+
\gamma \,  \varepsilon \,\displaystyle\int_{\partial\Omega}
\vert \nabla Q \vert^{p-2} \,
(\partial_k(\nabla Q)  \, n_k ) :  \nabla Q
\, d\sigma_{{\small \x}}=\displaystyle\sum_{i=1}^5 I_i.
\end{array}
\]
From \cite{GGRBRM}, by using the Sobolev's inequality:
\begin{equation}\label{gb}
\Vert \nabla Q \Vert_{L^{3p}(\Omega)}^{3p}  \le  \Vert |\nabla
Q |^{p/2}  \Vert_{L^6(\Omega)}^{6}  \le C \, \left(
\Vert \nabla Q \Vert_{L^p(\Omega)}^{p}
+ \Vert \nabla
\left( \vert \nabla Q \vert^{p/2} \right)
\Vert_{L^2(\Omega)}^{2} \right)^3,
\end{equation}
 the $I_1$-term is bounded as:
$$
\begin{array}{rcl}
I_1
% C_p \, \displaystyle\int_{\Omega} \vert \nabla {\u} \vert \vert \nabla Q \vert^p \, d\Omega
 & \le & C \, \Vert
\nabla {\u} \Vert_{L^q(\Omega)} \Vert \nabla Q
\Vert_{L^{pq/(q-1)}(\Omega)}^p \le  C \, \Vert \nabla {\u}
\Vert_{L^q(\Omega)} \Vert \nabla Q
\Vert_{L^p(\Omega)}^{\frac{p(2q-3)}{2q}} \Vert \nabla
Q \Vert_{L^{3p}(\Omega)}^{\frac{3p}{2q}} \\
\noalign{\vspace{-1ex}}\\
%& \le & \Vert \nabla {\u} \Vert_{L^q(\Omega)} \Vert \nabla {\bf
%d} \Vert_{L^p(\Omega)}^{p\alpha} \Vert
%\nabla (\vert \nabla Q \vert^{p/2}) \Vert_{L^2(\Omega)}^{2(1-\alpha)} \\
%\noalign{\vspace{-1ex}}\\
& \le & \displaystyle\frac{\gamma \, \varepsilon}{6} \left(\displaystyle\frac{2}{p}\right)^2 \, (p-2) \, \Vert \nabla (\vert \nabla Q
\vert^{p/2}) \Vert_{L^2(\Omega)}^2 +C_{\gamma,\varepsilon,p} \, \left( \Vert
\nabla {\u} \Vert_{L^q(\Omega)} + \Vert
\nabla {\u} \Vert_{L^q(\Omega)}^{\frac{2q}{2q-3}} \right) \, \Vert \nabla Q
\Vert_{L^p(\Omega)}^p
\end{array}
$$

\bigskip

By using (\ref{gb}) and   that $Q\in L^{\infty}(0,T;\mathbb{L}^{\infty}(\Omega))$ (given by a maximum principle), the $I_2$-term related to the stretching term $S(\nabla {\u},Q)$ can be bounded  as follows:

$$
\begin{array}{rcl}
I_2
& \le &   \displaystyle\int_{\Omega}
 \vert \nabla {\u} \vert
 \vert Q \vert
 \vert
\nabla Q \vert^{p-2} \vert D^2 Q \vert \, d{\bf
x} \\
\noalign{\vspace{-1ex}}\\
& \le & \Vert \nabla {\u} \Vert_{L^q(\Omega)} \left(
\displaystyle\int_{\Omega} \vert \nabla Q \vert^{p-2} \vert
D^2 Q \vert^2 \, d\x \right)^{1/2}
\Vert Q \Vert_{\mathbb{L}^{\infty}(\Omega)}
\Vert \nabla Q
\Vert_{L^{\frac{q(p-2)}{q-2}}(\Omega)}^{\frac{p-2}{2}}
\\
\noalign{\vspace{-1ex}}\\
& \le & C \, \Vert \nabla {\u} \Vert_{L^q(\Omega)} \left(
\displaystyle\int_{\Omega} \vert \nabla Q \vert^{p-2} \vert
D^2 Q \vert^2 \, d\x \right)^{1/2} \Vert \nabla Q
\Vert_{L^{p}(\Omega)}^{\frac{p(q-3)+q}{2q}} \Vert \nabla Q
\Vert_{L^{3p}(\Omega)}^{\frac{3(p-q)}{2q}}\\
\noalign{\vspace{-1ex}}\\
& \le & \displaystyle\frac{\gamma \,  \varepsilon \, (p-2)}{6} \, \left( \displaystyle\frac{2}{p}\right)^2 \Vert \nabla (\vert \nabla Q \vert^{p/2}) \Vert_{L^2(\Omega)}^2 \, + \displaystyle\frac{\gamma \, \varepsilon}{2} \, \displaystyle\int_{\Omega} \vert
\nabla Q \vert^{p-2} \vert D^2 Q \vert^2 \, d\x \\
\noalign{\vspace{-1ex}}\\
& + &  C_{\gamma \, \varepsilon,p} \,  \left\{\Vert \nabla Q
\Vert_{L^{p}(\Omega)}^p \, \left(\Vert \nabla {\u}
\Vert_{L^q(\Omega)}^{\frac{2q}{2q-3}} + 1
\right) +\Vert \nabla {\u}
\Vert_{L^q(\Omega)}^{\frac{q(p-1)}{2q-3}} \right\}.
\end{array}
$$
In order to know if additional regularity for $\nabla \u$ suffices to obtain  non-hilbertian estimates, we note that $L^{\frac{q(p-2)}{q-2}}(\Omega)$ interpolates between
$L^p(\Omega)$ and $L^{3p}(\Omega)$ if $q\le p \le
q/(3-q)$ (therefore $2\le q\le 3$). The $I_5$-term will be bounded below (see (\ref{nm1}) in Subsection \ref{bt}).

\bigskip

Collecting all the previous estimates, we arrive to the following inequality:
$$
\begin{array}{l}
\displaystyle\frac{2}{p}\displaystyle\frac{d}{dt} \Vert \nabla Q \Vert_{L^p(\Omega)}^p
  +  \gamma \, \varepsilon \displaystyle\int_{\Omega}
\vert D^2 Q \vert^2 \vert \nabla Q \vert^{p-2}
\, d\x +c \, \gamma \, \displaystyle\int_{\Omega} \left[
\vert Q \vert^2 \vert \nabla Q \vert^p + 2 \, \vert \nabla Q \vert^{p-2} \, (Q: \nabla Q)^2
\right] \, d\x\\
\noalign{\vspace{-1ex}}\\
\qquad +  \gamma   \varepsilon \, (p-2) \, \left( \displaystyle\frac{2}{p}\right)^2 \Vert \nabla (\vert \nabla Q \vert^{p/2}) \Vert_{L^2(\Omega)}^2\\
\noalign{\vspace{-1ex}}\\
\qquad \le    \widetilde{C} \, \left\{ \left( 1+ \Vert
\nabla {\u} \Vert_{L^q(\Omega)} + \Vert
\nabla {\u} \Vert_{L^q(\Omega)}^{2q/(2q-3)} 
 \right) \, \Vert \nabla Q
\Vert_{L^p(\Omega)}^p  + \Vert \nabla {\u}
\Vert_{L^q(\Omega)}^{q(p-1)/(2q-3)} \right\}.
\end{array}
$$
Taking into account that $q(p-1)/(2q-3) \le 2q/(2q-3)$ iff
$p \le 3$, the non-hilbertian estimates for $\nabla Q$ can be obtained for $p=3$, that is, $
\nabla Q \in {L}^{\infty}(0,T;\mathbb{L}^3(\Omega)) \cap
L^3(0,T;\mathbb{L}^9(\Omega))$.

\subsection{Special treatment for the boundary integral term $I_5$}\label{bt}
%In order to apply the Gronwall Lemma, we have to bound the boundary $I_5$-term of (\ref).
%\begin{equation}\label{frontera}
%\displaystyle\int_{\Gamma}
% \vert \nabla Q \vert^{p-2} \,
%\partial_{kl}^2 Q : \partial_k Q \, n_l
%\, d\sigma
%\end{equation}
By using a boundary parametrization, we study the case where $\partial\Omega$ is a finite union of $\Gamma_i$, $\partial\Omega=\sum_{i=1}^I\Gamma_i$, each one defined through a parametrization of one variable of $\mathbb{R}^3$ in function of the other two. For example, choosing a parametrization where the  $z$ variable is function of $\x=(x,y)$, we have
$\Gamma_i=\{(\x,z) \in \mathbb{R}^3, \, z=f(\x), \x=(x,y)\in \omega\}$. Denoting $\partial_1=\partial_{x}$, $\partial_2=\partial_{y}$, $\partial_3=\partial_{z}$, the normal vector ${\bf n}=(n_1,n_2,n_3)$ to $\Gamma_i$ is proportional to ${\bf m}=(m_1,m_2,-1)=(\partial_1 f(\x),\partial_2 f(\x),-1)$. Observe that in this case, 
$$
\partial_{\bf n} Q \vert_{\partial\Omega}=0\quad \mbox{iff} \quad -\partial_3 Q (\x,f(\x))+
\partial_1 Q(\x,f(\x)) \, m_1+ \partial_2 Q(\x,f(\x)) \, m_2={\bf 0}
$$ 
for $\x\in \omega
$. Hence
\begin{equation}\label{dq0}
\partial_3 Q = \partial_1 Q \, m_1 + \partial_2 Q \, m_2  , \quad \x\in \omega
\end{equation}
In this boundary $\Gamma_i$, we can derivate with respect to $x$ and $y$, obtaining (denoting $m_{j,i}=\partial_i m_j$, $j=1,2$; $i=1,2,3$):
\begin{eqnarray}
\partial_{31}^2 Q + \partial_{33}^2 Q \, m_1  =
(\partial_{11}^2 Q + \partial_{13}^2 Q \, m_1) \, m_1
+ \partial_1 Q \, m_{1,1} +  (\partial_{21}^2 Q + \partial_{23}^2 Q \, m_1) \, m_2
+ \partial_2 Q \, m_{2,1} \label{dq1} \\
%\noalign{\vspace{-2ex}} \nonumber \\
\partial_{32}^2 Q + \partial_{33}^2 Q \, m_2  =
(\partial_{12}^2 Q + \partial_{13}^2 Q \, m_2) \, m_1
+ \partial_1 Q \, m_{2,1} +  (\partial_{22}^2 Q + \partial_{23}^2 Q \, m_2) \, m_2
+ \partial_2 Q \, m_{2,2}\label{dq2}
\end{eqnarray}
By multiplying (\ref{dq1}) by $ \partial_1 Q$ and  (\ref{dq2}) by $\partial_2 Q$ and using (\ref{dq0}), we obtain for the left hand-side:
$$
\partial_{31}^2 Q : \partial_1 Q +
\partial_{33}^2 Q :  \partial_1 Q \, m_1 +
 \partial_{32}^2 Q : \partial_2 Q +
 \partial_{33}^2 Q :  \partial_2 Q \, m_2  =
 \displaystyle\sum_{k=1}^3 \, \partial_{3k}^2 Q : \partial_k Q
$$
and for the right-hand side (using (\ref{dq0}) and the fact that $m_{1,2}=m_{2,1}$):
$$
\begin{array}{l}
 \displaystyle\sum_{k=1}^2 \, \partial_{1k}^2 Q : \partial_k Q \, m_1 +
  \partial_{13}^2 Q :
  \left(\partial_1 Q \, m_1+ \partial_2 Q \, m_2 \right) \, m_1\\
  \noalign{\vspace{-2ex}}\\
  \qquad +  \displaystyle\sum_{k=1}^2 \, \partial_{2k}^2 Q : \partial_k Q \, m_2 +
  \partial_{23}^2 Q :
  \left( \partial_1 Q \, m_1+ \partial_2 Q \, m_2 \right) m_2\\
  \noalign{\vspace{-2ex}}\\
\qquad   + \vert \partial_1 Q \vert^2 \, m_{1,1} +
  2 \, \partial_1 Q : \partial_2 Q \, m_{1,2} +
  \vert \partial_2 Q \vert^2 \, m_{2,2}\\
 \noalign{\vspace{-2ex}}\\
\qquad =  \displaystyle\sum_{k=1}^3 \, \partial_{1k}^2 Q : 
\partial_k Q \, m_1 + \displaystyle\sum_{k=1}^3 \,
\partial_{2k}^2 Q : \partial_k Q \, m_2
  + \vert \partial_1 Q \vert^2 \, m_{1,1} +
  2 \, \partial_1 Q : \partial_2 Q \, m_{1,2} +
  \vert \partial_2 Q \vert^2 \, m_{2,2}.
\end{array}
$$
Then,
\begin{equation}\label{estrella2-1}
\begin{array}{rcl}
0 & =  & \displaystyle\sum_{k=1}^3 \, \partial_{3k}^2 Q
: \partial_k Q  \, n_3 +
 \displaystyle\sum_{k=1}^3 \, \partial_{1k}^2 Q : \partial_k Q \, n_1 +
\displaystyle\sum_{k=1}^3 \, \partial_{2k}^2 Q :
\partial_k Q \, n_2\\
\noalign{\vspace{-1ex}}\\
& + &
 \vert \partial_1 Q \vert^2 \, \displaystyle\frac{m_{1,1}}{\vert {\bf m} \vert}
 + 2 \, \partial_1 Q : \partial_2 Q \, \displaystyle\frac{m_{1,2}}{\vert {\bf m} \vert}
 + \vert \partial_2 Q \vert^2 \, \displaystyle\frac{m_{2,2}}{\vert {\bf m} \vert}.
\end{array}
\end{equation}

\bigskip

Therefore, using (\ref{estrella2-1}) in the boundary $I_5$ term, we have:
\begin{equation}\label{nm1}
\begin{array}{rcl}
I_5 & = &
 \gamma   \varepsilon  \displaystyle\int_{\partial\Omega}
 \vert \nabla Q \vert^{p-2} 
\partial_{kl}^2 Q : \partial_k Q \, n_l
\, d\sigma  \\
\noalign{\vspace{-2ex}}\\
& = &  \gamma   \varepsilon  \displaystyle\int_{\partial\Omega}
\vert \nabla Q \vert^{p-2}  \left( \partial_{k3}^2
Q :  \partial_k Q \, n_3 +  \partial_{k1}^2
Q :  \partial_k Q \, n_1 +  \partial_{k2}^2
Q :  \partial_k Q
\, n_2 \right) \, d\sigma\\
\noalign{\vspace{-2ex}}\\
& = & - \gamma  \,  \varepsilon \, \displaystyle\int_{\partial\Omega}
\vert \nabla Q \vert^{p-2} \, \left( \vert \partial_1
Q \vert^2 \, m_{1,1} + 2 \, \partial_1 Q :
\partial_2 Q \, m_{1,2} + \vert \partial_2 Q
\vert^2 \, m_{2,2}
 \right)/\vert {\bf m} \vert \, d\sigma \\
 \noalign{\vspace{-2ex}}\\
 & \le &
 C \, \gamma \,  \varepsilon \, \displaystyle\int_{\partial\Omega} \vert \nabla Q \vert^p \, d\sigma =
 C \, \gamma \,  \varepsilon \, \Vert \nabla Q \Vert_{L^p(\partial\Omega)}^p
 = C \, \gamma \,  \varepsilon \, \Vert \vert \nabla Q \vert^{p/2} \Vert_{L^2(\partial\Omega)}^2
 \le  \widetilde{C} \, \gamma \,  \varepsilon \, \Vert \vert \nabla Q \vert^{p/2}
 \Vert_{H^{1/2+\epsilon}(\Omega)}^2 \\
 \noalign{\vspace{-2ex}}\\
 & \le & \gamma \,  \varepsilon \, \delta \Vert \nabla \left( \vert \nabla Q \vert^{p/2} \right) \Vert_{L^2(\Omega)}^2 + \gamma \,  \varepsilon \, C_{\delta} \Vert \vert \nabla Q \vert^{p/2} \Vert_{L^2(\Omega)}^2  \\
 \noalign{\vspace{-2ex}}\\
 & \le & \displaystyle\frac{\gamma\,  \varepsilon}{6}  \,
 \left(\displaystyle\frac{2}{p}\right)^2 \, (p-2) \,
 \Vert \nabla \left( \vert \nabla Q \vert^{p/2} \right) \Vert_{L^2(\Omega)}^2 +
 C_{p, \gamma,  \varepsilon} \, \Vert \nabla Q  \Vert_{L^p(\Omega)}^p
\end{array}
\end{equation}

\begin{remark} In the special case of  $\partial\Omega=\{
(\x,0) \in \mathbb{R}^3\}$, the normal vector is ${\bf n}=(n_1,n_2,n_3)=(0,0,-1)$ and
$
\partial_{\bf n} Q \vert_{\partial\Omega}= \, \partial_3  Q (\x,0)={\bf 0}$ which implies that $\partial_{3i}^2 Q (\x,0)={\bf 0}$ for $i=1,2$. Therefore, the boundary $I_5$ term in (\ref{nm1}) is given by $$
\displaystyle\int_{\partial\Omega}
\vert \nabla Q \vert^{p-2} \,
\partial_{kl}^2 Q \cdot  \partial_k  Q \, n_l
\, d\sigma= -\displaystyle\int_{\partial\Omega}
\vert \nabla Q \vert^{p-2} \,
\partial_{k3}^2 Q \cdot \partial_k Q
\, d\sigma =0 .$$
\end{remark}

\bigskip

\begin{remark}[Open problem for Dirichlet b.c.]
The extension to Dirichlet boundary conditions for $Q$ is an open problem. In fact, even in the case of $\partial\Omega=\{
(\x,0) \in \mathbb{R}^3\}$, since ${\bf n}=(0,0,-1)$, the boundary integral $I_5$ reduces to:
$$
\displaystyle\int_{\partial\Omega}
 \vert \nabla Q \vert^{p-2} \,
\partial_{kl}^2 Q \cdot \partial_k Q \, n_l
\, d\sigma=-\displaystyle\int_{\partial\Omega}
\vert \nabla Q \vert^{p-2} \,
\partial_{33}^2 Q \cdot  \partial_3 Q
\, d\sigma =-\displaystyle\frac{1}{2} \, \displaystyle\int_{\partial\Omega}
\vert \nabla Q \vert^{p-2} \, \partial_3 \left( \vert \partial_3 Q \vert^2 \right)
\, d\sigma.
$$
But,  we do not have any information for $\partial_3 Q$ on $\partial\Omega$, hence
we do not know how to bound this boundary term in function of $\Vert \nabla Q \Vert_{L^p(\Omega)}$ and $\Vert \nabla Q \Vert_{L^{3p}(\Omega)}$.
\end{remark}

\section{Proof of Theorem \ref{t2}.}\label{sec-reg}
Regularity criteria for $\nabla \u$ given  in (\ref{ar-u}) and (\ref{reg-ad-u}), respectively, joint to $Q \in L^{\infty}((0,T) \times\Omega )$ given by the maximum principle for $Q$,  imply that  $
S(\nabla \u, Q) \in L^{2q/(2q-3)}(0,T;{\bf L}^{q}(\Omega))$ for $2 \le q\le 3$ or $3/2 \le q\le 3$, respectively.
In the same way, as $\nabla Q   \in  L^{\infty}(0,T;\mathbb{L}^3(\Omega))$ and $\u  \in  L^{2q/(2q-3)}(0,T;{\bf W}^{1,q}(\Omega))$, then the convective term $
(\u \cdot \nabla) Q \in L^{2q/(2q-3)}(0,T;\mathbb{L}^{q}(\Omega))$.
It is easy to prove that $f(Q) \! \in \! L^{2q/(2q-3)}(0,T;\mathbb{L}^{q}(\Omega))$. Indeed, the worse term in (\ref{ff}) is bounded by:
$$
\begin{array}{rcl}
\Vert \vert Q \vert^3 \Vert_{\mathbb{L}^{q}(\Omega)}^{2q/(2q-3)} 
& \le & \Vert Q \Vert_{\mathbb{L}^{3q}(\Omega)}^{6q/(2q-3)}
\le C_0 \, \Vert Q \Vert_{\mathbb{W}^{1,3q/(q+1)}(\Omega)}^{ 6q/(2q-3) }\\
\noalign{\vspace{-1ex}}\\
& \le & 
\left\{\begin{array}{ll}
C_0 \, \Vert Q \Vert_{\mathbb{H}^{1}(\Omega)}^{6q/(2q-3)} & \mbox{if $q\le 2$,}
\\
\noalign{\vspace{-2ex}}\\
C_0 \, \Vert Q \Vert_{\mathbb{H}^1(\Omega)}^{3(q+2)/(2q-3)}
\Vert Q \Vert_{\mathbb{H}^2(\Omega)}^{3(q-2)/(2q-3)} & \mbox{if $q>2$,}
\end{array}\right.
\end{array}
$$
where $C_0$ is the constant of the embedding $\mathbb{W}^{1,3q/(q+1)}(\Omega) \hookrightarrow \mathbb{L}^{3q}(\Omega)$. Since $Q \in L^{\infty}(0,T;\mathbb{H}^1(\Omega)) \cap L^2(0,T;\mathbb{H}^2(\Omega))$ and $3(q-2)/(2q-3)\le 2$, therefore $f(Q) \in L^{2q/(2q-3)}(0,T;\mathbb{L}^{q}(\Omega))$. In summary,
\begin{equation}\label{mars-1}
\mathbb{G}=S(\nabla \u,Q)-(\u \cdot \nabla) Q -f(Q) \in L^{2q/(2q-3)}(0,T;\mathbb{L}^{q}(\Omega)).
\end{equation}
Depending on the value of $q$,  two situations can be presented:
\begin{eqnarray}
 \mbox{for $3/2\le q\le 5/2$,} \quad \mathbb{G} \in L^\gamma(0,T;\mathbb{L}^\gamma(\Omega))  \quad  \mbox{where } \gamma=q=\min \left\{q,\frac{2q}{2q-3}\right\}\in [3/2,5/2] ,\label{mars-2} \\
\nonumber \\
\mbox{for  $5/2\le q\le 3$,} \quad   \mathbb{G} \in L^\gamma(0,T;\mathbb{L}^\gamma(\Omega)) \quad  \mbox{where } \gamma=\frac{2q}{2q-3}=
\min \left\{q,\frac{2q}{2q-3}\right\}\in [2,5/2] .\label{mars-3}
\end{eqnarray}
Now, we study the regularity of the parabolic problem given by the  heat equation under homogeneous Neumann boundary conditions:
\begin{equation}\label{2bvp}
\partial_t w - \gamma \, \varepsilon \, \Delta w = f \quad
\mbox{in $\Omega \times (0,T)$,} \qquad
\partial_{\bf n} w \vert_{\Omega}={ 0} \quad   \mbox{in $(0,T)$,} \qquad w\vert_{t=0}=w_0 \, \, \mbox{in $\Omega$.}
\end{equation}
 We rewrite Theorem 17 of the  Solonnikov's paper \cite{Sol64} for the case of Neumann boundary conditions. 
\begin{theorem}[Theorem~17 of \cite{Sol64}] \label{t-3}
Let $\Omega \subset \mathbb{R}^3$ be a bounded domain with boundary $\partial\Omega$ of class $C^{2+\epsilon}$ for some $\epsilon>0$. Let $\gamma\in (1,+\infty)$ and assume that $w_0 \in B_\gamma^{2-2/\gamma}(\Omega)$ and $f\in L^\gamma(0,T;L^\gamma(\Omega))$. Then the solution $w$ of (\ref{2bvp}) satisfies:
$$
%\begin{equation}\label{2bvp-e1}
\Vert w \Vert_{L^\gamma(0,T;W^{2,\gamma}(\Omega))} + \Vert \partial_t w \Vert_{L^\gamma(0,T;L^\gamma(\Omega))} \le C \, \left(
\Vert f \Vert_{L^\gamma(0,T;L^\gamma(\Omega))} + \Vert w_0 \Vert_{B_\gamma^{2-2/\gamma}(\Omega)}
\right).
$$
%\end{equation}
\end{theorem}
Now, we apply Theorem \ref{t-3} to any component $Q_{ij}$ of $Q$,  which is the solution of the problem (\ref{2bvp}) for $f= \mathbb{G}_{ij}$ and  $w=(Q_0)_{ij}$,
where $\mathbb{G}$ is defined in (\ref{mars-1}). According to (\ref{mars-2}) and (\ref{mars-3}), then $Q \in L^{\gamma}(0,T;\mathbb{W}^{2,\gamma}(\Omega))$ and 
 $\partial_t Q \in L^{\gamma}(0,T;\mathbb{L}^{\gamma}(\Omega))$
 with $\gamma$ given in (\ref{gamma}), and the proof of  Theorem~\ref{t2} is finished.

\section{Extension to the traceless and/or symmetric models}

In this section, we analyzes the influence on the proofs of Theorems \ref{t1} and \ref{t2} of both constraints, either symmetry (replacing the stretching term $S(\nabla \u,Q)$ defined in (\ref{def-S}) by (\ref{S-2})), or tracelessness (replacing the function $f(Q)$ defined in (\ref{ff}) by (\ref{fz})). The proofs of Corollaries \ref{c1} and \ref{c2} do not need be modified.

\subsection{Effect of the symmetry.}

When the definition (\ref{S-2}) for $S(\nabla \u,Q)$ is considered, the proof of Theorem \ref{t1} must be modified only  replacing $\nabla\u$ by $\W$ in the $I_2$-term depending on the stretching.

The proof of Theorem \ref{t2} must be revised when the second-member $\mathbb{G}$ is defined in (\ref{mars-1}). Now, 
$$
\mathbb{G}=S(\W,Q)-(\u \cdot \nabla) Q -f(Q) .
$$
Then, hypothesis (\ref{ar-u}) also implies that $\mathbb{G}  \in L^{2q/(2q-3)}(0,T;\mathbb{L}^{q}(\Omega))$, hence the proof of Theorem \ref{t2} is deduced.
\subsection{Effect of the tracelessness.}
Instead of seeing the effect of the new definition of $f(Q)$ given in (\ref{fz}) over Theorems \ref{t1} and \ref{t2}, we will study the influence of replacing (\ref{ff}) by the more general case:
\begin{equation}\label{m-1}
f(Q)= a \, Q - \displaystyle\frac{b}{3} \, \left(
Q^2+QQ^t+Q^tQ \right )+c  \, \vert Q \vert^2 \, Q -\alpha(Q) \, \mathbb{I}\quad
\mbox{ with  $a$, $b\in \mathbb{R}$ and $c
>0$}
\end{equation}
with $\alpha(Q)$ a convex combination of 
$
\alpha_1(Q):=\frac13\left(-a \,
tr(Q)+\frac{b}{3} \,  \left(
Q:Q^t+ 2 \vert Q \vert^2 \right) \right)
$ and $
\alpha_2(Q):=-\frac{tr(f(Q))}{3}$. Observe that, since $tr(Q)=0$ can be deduced (\cite{GGRBinprep2}), then $\alpha(Q)$ reduces to:
\begin{equation}\label{alfa-3}
\alpha(Q):=\frac{b}{9} \,  \left(
Q:Q^t+ 2 \vert Q \vert^2 \right) .
\end{equation}

In order to see the effect of \eqref{m-1}-\eqref{alfa-3} on the proof of Theorem \ref{t1}, it suffices to analyze the new term appearing when $ -\nabla \cdot \left(
\vert \nabla Q \vert^{p-2} \nabla Q \right)$ is taken as test function in the $Q$-system of (\ref{modelo-u-1}):
$$
\begin{array}{rcl}
\displaystyle\int_{\Omega}
\alpha(Q) \mathbb{I} : \left[ \nabla \cdot \left(
\vert \nabla Q \vert^{p-2} \nabla Q \right) \right]
\, d\x & = & -\displaystyle\int_{\Omega}
\nabla (\alpha(Q)) \, \vert \nabla Q \vert^{p-2} \nabla (tr(Q)) 
\, d\x \\
\noalign{\vspace{-1ex}}\\
& + & \displaystyle\int_{\partial\Omega}
\alpha(Q) \, \vert \nabla Q \vert^{p-2} \partial_{\bf n} (tr(Q))
\, d\sigma=0
\end{array}$$
and therefore, no modifications in the proof of Theorem \ref{t1} must be done.

Concerning the proof of Theorem \ref{t2}, again we have to look at the term $\mathbb{G}=S(\nabla \u,Q)-(\u \cdot \nabla) Q -f(Q)$ with $f(Q)$ given by (\ref{m-1}) with $\alpha(Q)$ defined in (\ref{alfa-3}). Observe that:
$$
\vert \alpha(Q) \vert \le \displaystyle\frac{b}{3} \, \vert Q \vert^2
$$
and therefore for any $q\le 3$:
$$
\Vert \vert Q \vert^2 \Vert_{\mathbb{L}^q(\Omega)}^{2q/(2q-3)} =\Vert Q \Vert_{\mathbb{L}^{2q}(\Omega)}^{4q/(2q-3)}
\le
C_0 \, \Vert Q \Vert_{\mathbb{H}^1(\Omega)}^{4q/(2q-3)}.
$$
Since $Q \in L^{\infty}(0,T;\mathbb{H}^1(\Omega))$, then $\alpha(Q) \, \mathbb{I}\in L^{\infty}(0,T;\mathbb{L}^q(\Omega))$ and therefore $\mathbb{G} \in L^{2q/(2q-3)}(0,T;\mathbb{L}^q(\Omega))$. From this point, the proof of Theorem \ref{t2} does not change.



%% References


\bibliographystyle{model1-num-names}
\bibliography{<your-bib-database>}



\end{document}